\documentclass[11pt]{article}
\usepackage{amssymb,latexsym}
\usepackage{epsfig}
\usepackage{eufrak}
\usepackage{amsmath}
\usepackage{mathrsfs}
\usepackage{color}

\newtheorem{theorem}{Theorem}
\newtheorem{lemma}{Lemma}

\newcommand{\be}{\begin{equation}}
\newcommand{\ee}{\end{equation}}
\newcommand{\bee}{\begin{eqnarray*}}
\newcommand{\eee}{\end{eqnarray*}}
\newcommand{\bel}{\begin{eqnarray}}
\newcommand{\eel}{\end{eqnarray}}
\newcommand{\bec}{\begin{cases}}
\newcommand{\eec}{\end{cases}}
\newcommand{\bem}{\begin{bmatrix}}
\newcommand{\eem}{\end{bmatrix}}

\newcommand{\bed}{\begin{description}}
\newcommand{\eed}{\end{description}}
\newcommand{\bei}{\begin{itemize}}
\newcommand{\eei}{\end{itemize}}
\newcommand{\ben}{\begin{enumerate}}
\newcommand{\een}{\end{enumerate}}

\newcommand{\beL}{\begin{lemma}}
\newcommand{\eeL}{\end{lemma}}
\newcommand{\beT}{\begin{theorem}}
\newcommand{\eeT}{\end{theorem}}

\newcommand{\bpf}{\begin{pf}}
\newcommand{\epf}{\end{pf}}

\newcommand{\pfbox}{\hfill\mbox{$\Box$}}

\newenvironment{pf}{\paragraph*{Proof{\rm.}}}{\pfbox\bigskip}

\begin{document}

\title{{\bf Fast Construction of Robustness Degradation Function \thanks{Received by the editors
June 5, 2002; accepted for publication (in revised form) August 6,
2003; published electronically January 22, 2004. This research was
supported in part by grants from NASA (NCC5-573) and LEQSF
(NASA/LEQSF(2001-04)-01).  }}}
\author{Xinjia Chen, Kemin Zhou and Jorge L. Aravena\\
Department of Electrical and Computer Engineering\\
Louisiana State University\\
Baton Rouge, LA 70803}

\date{Received in June 5, 2002;  Revised in August 6, 2003}

\maketitle

\begin{abstract}  We develop a fast algorithm to
construct the robustness degradation function, which describes quantitatively
the relationship between
the proportion of systems
guaranteeing  the robustness requirement and the radius of the uncertainty set.
This function can be applied to predict whether a
controller design based on an inexact mathematical
model will perform satisfactorily when implemented on the true system.

\end{abstract}

\section{Introduction}

In recent years, there has been growing interest
in the development of probabilistic methods for robustness analysis and design
problems aimed at overcoming the computational complexity
and the conservatism issue of the deterministic
worst-case framework \cite{RS, SR, MS, KT, TD, BLT, BL, BP, Polyak, CDT1, CDT2, C1, C2,
vind, vind2, SB, BLT1}.  In the deterministic worst-case framework,
one is interested in knowing if the robustness requirement is
guaranteed for every value of the uncertainty.
However, it should be borne in mind that the uncertainty set
may include worst cases which never happen in reality.
Instead of seeking the worst-case guarantee, it is sometimes ``acceptable''
 that the robustness requirement is satisfied for most of the cases.
 It has been demonstrated that the proportion of systems
 guaranteeing  the robustness requirement
 can be close to $1$ even if the radii of the uncertainty set are
 much larger than the worst-case deterministic robustness margin \cite{ BLT, BLT2, BP, CDT1, Polyak}.
 Therefore, it is of practical importance to construct a function
 which describes quantitatively the relationship between
the proportion of systems
guaranteeing  the robustness requirement and the radius of the uncertainty set.
Such a function can serve as a guide for control engineers in
evaluating the robustness of
a control system once a controller design is completed.
Such a function, referred to as a {\it robustness degradation function},  has been proposed
by a number of researchers \cite{BLT, CDT1}.
For example,  Barmish, Lagoa, and Tempo \cite{BLT} have
constructed a curve of the
robustness margin amplification versus risk in a probabilistic setting.
In a similar spirit, Calafiore, Dabbene, and Tempo \cite{CDT1, CDT2} have
constructed a probability degradation function in the context of
 real and complex parametric uncertainty.

In this paper, allowing the robustness analysis to be performed in a distribution-free manner,
we introduce the concept of {\it proportion} and adopt the assumption from the classical robust control framework that
 uncertainty is deterministic and bounded.
It follows naturally that the robustness of a system can be
reasonably measured by the ratio of the volume (Lebesgue measure) of the set of
uncertainty guaranteeing the robustness requirement
to the overall set of uncertainty \cite{TD}.  Evaluation of such a measure of robustness requires
generating samples with uniform distribution over uncertainty sets such as a
spectral normal ball or  an $l_p$ ball.  The difficulty
 of generating such samples has been successfully resolved in \cite{CDT1, CDT2}.

The conventional method for constructing the robustness function is to perform,
independently,  a certain number of
simulations for each value of the uncertainty radius and then plot the function.
Although such a curve can be applied to evaluate the robustness of the control system,
it may be computationally expensive.
This is especially true when many cycles of controller synthesis and robustness analysis are
needed in the development of a high performance control system.
Motivated by this situation, we focus on the machinery that
can make the construction of such a function efficient.
  We have developed a sample reuse algorithm that allows
the simulations to be conducted in an iterative manner.
The idea is to start simulation from the larger uncertainty set and save appropriate evaluations of the robust
requirement for the use of later simulations on the smaller uncertainty set.
In this way the total number of simulations
can be reduced significantly as compared to the conventional method.

In addition to deriving our sample reuse algorithm from the worst-case deterministic framework,
we show that the technique is also applicable when considering the random nature of the uncertainty.
In such cases, the worst-case properties of uniform distribution given in the
pioneering work \cite{BL, BLT, bai} allow our algorithm
to be applied to efficiently solve a wide variety of robustness analysis problems.
In particular, the radial truncation theory \cite{BLT}
can be applied to robustness analysis
problems with uncertainty bounding sets defined as spectral norm balls and $l_p$ balls.

The organization of the paper is as follows.  Section 2 gives the problem formulation.
Section 3 presents our sample reuse algorithm. Section 4 is the performance analysis of the algorithm.
Section 5 applies the algorithm to examples.
Section 6 shows the justification
of the algorithm for the case of random uncertainties.
Section 7 is the conclusion.  The proofs of the theorems are included as an appendix.

\section{Problem formulation}

We adopt the assumption, from the classical robust control framework, that the
uncertainty is deterministic and bounded.
We formulate a general robustness analysis problem as follows.

Let ${\bf P}$ denote a robustness requirement.
The definition of ${\bf P}$ can be a fairly complicated
combination of the following:

\begin{itemize}

\item stability or ${\cal D}$-stability;

\item $H_{\infty}$ norm of the closed-loop transfer function;

\item  time specifications such as
overshoot, rise time, settling time, and steady state error.

\end{itemize}

Let ${\cal B}(r)$ denote the
set of uncertainties with size smaller than $r$. In applications, we are usually dealing with
uncertainty sets such as the following:
\begin{itemize}
\item $l_p$ ball
$
{\cal B}_p(r) :=\left\{ \Delta \in {\bf R}^n:\;||\Delta||_p \leq r \right\},
$
where $||.|_p$ denotes the $l_p$ norm and $p = 1, 2, \ldots,\infty$.
In particular, ${\cal B}_{\infty} (r)$ denotes a box.
\item  Spectral norm ball
$
{\cal B}_{\sigma}(r) := \{ \Delta \in {\bf \Delta}: \bar{\sigma}(\Delta) \leq r \},$
where $\bar{\sigma}(\Delta)$ denotes the largest singular value of $\Delta$.  The
 class of allowable perturbations is \vspace{3pt}
\begin{equation}
\label{str}
{\bf \Delta} : =\{ {\rm blockdiag} [q_1I_{r_1}, \ldots, q_sI_{r_s}, \Delta_1, \ldots, \Delta_c  ]\},
\end{equation}
\vspace{3pt} where $q_i \in \Bbb{F}, \;i=1, \ldots,s$ are scalar parameters with multiplicity $r_1,\ldots,r_s$ and
$\Delta_i \in \Bbb{F}^{n_i \times m_i}, \;i = 1, \ldots, c$ are possibly repeated full blocks.
Here $\Bbb{F}$ is either the complex field ${\bf C}$ or the real field ${\bf R}$.

\item Homogeneous star-shaped bounding set
$
{\cal B}_H (r):=\left\{ r (\Delta - \Delta_0) + \Delta_0 :
 \Delta \in Q \right\},
 $
 where $Q \subset {\bf R}^n$ and $\Delta_0 \in Q$
 (see \cite{BLT} for a detailed illustration).
\end{itemize}
Throughout this paper, ${\cal B} (r)$ refers to any type of uncertainty set described above.
Define a function $\ell ( . )$  such that,
 for any $X$,\vspace{3pt}
 \[
 \ell ( X ) := \min \{ r : X \in {\cal B} (r) \},
 \]
 \vspace{3pt} i.e.,
${\cal B}(\ell ( X ))$ includes $X$ exactly in the boundary.  By such definition,\vspace{3pt}
\[
\ell ( X ) = \min \left\{ r :  \frac{X- \Delta_0}{r} + \Delta_0 \in Q \right\},
\]\vspace{3pt}
\[
\ell ( X ) = \bar{\sigma}(X),
\]\vspace{3pt}
and\vspace{3pt}
\[
\ell ( X ) = ||X||_p
\]\vspace{3pt}
in the context of a homogeneous star-shaped bounding set, spectral norm ball, and $l_p$ ball, respectively.

To allow the robustness analysis to be performed in a distribution-free manner,
we introduce the notion of {\it proportion} as follows.  For any $\Delta \in {\cal B}(r)$
there is an associated system $G(\Delta)$. We define {\it proportion} as follows:\vspace{3pt}
\[
\Bbb{P} (r) := \frac{{\rm vol} ( \{\Delta \in {\cal B}(r):
{\rm The \; associated \; system}
\; G(\Delta) \; {\rm guarantees} \; {\bf P} \}   ) } { {\rm vol} ( {\cal B}(r) ) }
\]\vspace{3pt}
with \vspace{3pt}
\[
{\rm vol} (S) := \int_{q \in S} dq,
\]\vspace{3pt}
where the notion of $dq$ is illustrated as follows:

\begin{itemize}

\item (I): If $q =[x_{rs}]_{n \times m}$ is a real matrix in ${\bf R}^{n \times m}$,
then $dq = \prod_{r=1}^n \prod_{s=1}^m dx_{rs}$.

\item (II): If $q =[x_{rs} + j y_{rs}]_{n \times m}$ is a complex
 matrix in ${\bf C}^{n \times m}$, then $dq = \prod_{r=1}^n \prod_{s=1}^m ( dx_{rs} dy_{rs} )$.

\item (III): If $q \in {\bf \Delta}$, i.e., $q$ possesses a block structure defined by (\ref{str}), then
$dq = (\prod_{i=1}^s d q_i) (\prod _{i=1}^c d \Delta_i)$, where the notion
of $dq_i$ and $d \Delta_i$ is defined by (I) and (II).

\end{itemize}

It follows that $\Bbb{P} (r)$ is a reasonable measure of the robustness of the system \cite{CDT1,TD2}.
In the worst-case deterministic framework,
we are interested only in knowing if ${\bf P}$ is guaranteed for every $\Delta$.
However, one should bear in mind that the uncertainty set in
our model may include worst cases which never happen in reality.
Thus, it would be ``acceptable'' in many applications
if the robustness requirement ${\bf P}$ is satisfied for most of the cases.
Hence,
 due to the inaccuracy of the model, we should also
 obtain the value of $\Bbb{P} (r)$ for uncertainty radius $r$ which exceeds the deterministic robustness margin.

Clearly, $\Bbb{P} (r)$ is deterministic
in nature.  However, we can resort to
a probabilistic approach to evaluate $\Bbb{P} (r)$.  To see this, one needs to
observe that a random variable with
{\it uniform distribution} over ${\cal B} (r)$, denoted by ${\bf \Delta}^u$, guarantees that\vspace{3pt}
\[
\Pr \{ {\bf \Delta}^u \in S \} = \frac{{\rm vol} (S \bigcap {\cal B}(r)) } { {\rm vol} ({\cal B}(r))  }
\]
\vspace{3pt} for any $S$, and thus\vspace{3pt}
\[
\Bbb{P} (r) =  \Pr \{
{\rm The \; associated \; system}
\; G({\bf \Delta}^u) \; {\rm guarantees} \; {\bf P} \}.
\]
\vspace{3pt}  It follows that a Monte Carlo method can be employed to estimate $\Bbb{P} (r)$
based on independently and identically distributed (i.i.d.) observations of ${\bf \Delta}^u$.

It is interesting to know how the function $\Bbb{P} (r)$ degrades
with respect to $r$ when $r$ increases from $a$ to $b$, where $b, \;
a \geq 0$. In a similar spirit, such a function has been proposed as
a {\it confidence degradation function} in \cite{BLT} and as a {\it
probability degradation function} in \cite{CDT1, CDT2}. In this
paper, we refer to the function $\Bbb{P}(.)$ as a {\it robustness
degradation function} for the following reasons. First, we introduce
the confidence interval for assessing the accuracy of the estimate
of $\Bbb{P}(r)$. To be useful, every numerical method should be
associated with an assessment for the accuracy of the estimate.
Monte Carlo simulation is no exception. To avoid confusion, we
reserve the notion of ``confidence'' for the purpose of interval
estimation. Second, we introduce the concept of {\it proportion} for
measuring robustness, which has no probabilistic content. Third,
$\Bbb{P}(r)$ is a robustness measure and is usually decreasing
 with respect to $r$ when $\Bbb{P}(r)$ is close to $1$.

To construct such a function of practical importance, the
conventional way is to grid the interval $[a, b]$ as $a = \rho_1 <
\rho_2 < \cdots < \rho_l = b$ and estimate $\Bbb{P} (\rho_i)$ by
conducting $N$ i.i.d. sampling experiments for each $\rho_i$. In
total, we need $Nl$ samples.  In the next section we show that the
number of experiments can be significantly reduced.

\section{Sample reuse algorithm}

To improve efficiency, we shall make use of the following simple
yet important observation.

{\it Let $q^*$ be an observation of a random variable with
uniform distribution over ${\cal B}(r^*)  \supseteq {\cal B}(r)$
such that $q^* \in {\cal B}(r)$.  Then $q^*$ can also be viewed
as an observation of a random variable with uniform distribution
over ${\cal B}(r)$. }


In our algorithm, we flip the order of $\rho_i$ by defining\vspace{3pt}
\[
r_i = \rho_{l+1-i}
\]
\vspace{3pt} for $i = 1,2, \ldots,l$.  Thus, the direction of simulation is backward.
Our algorithm is described as follows.
\unskip

\underline{{\sc Sample Reuse Algorithm}}

\begin{itemize}
\item Input: Sample size $N$, confidence parameter $\delta \in (0,1)$
and uncertainty radii $r_i, \;\; i=1, 2, \ldots, l$.

\item Output: Proportion estimate $\widehat{\Bbb{P}}_i$ and the related confidence interval
for $i=1,\ldots,l$.  In the following, $m_{i1}$ denotes the number of sampling experiments
 conducted at $r_i$, and
 $m_{i2}$ denotes the number of observations guaranteeing ${\bf P}$
 during the $m_{i1}$ sampling experiments.

\item {{\it Step} $1$ (initialization).}
 Let $M =[m_{ij}]_{l \times 2}$ be a zero matrix.

\item {{\it Step} $2$ (backward iteration).} For $i=1$ to $i=l$ do the following:

      \begin{itemize}

      \item Let $r  \leftarrow  r_i$.

      \item While $m_{i1} < N$ do the following:

              \begin{itemize}

                \item Generate uniform
                  sample $q$ from ${\cal B}(r)$. Evaluate the robustness requirement
                  {\bf P} for $q$.

                \item Let $m_{s1} \leftarrow m_{s1}+1$ for any $s$ such that $ r \geq r_s \geq \ell ( q )$.

                \item If robustness requirement {\bf P} is satisfied
                      for $q$, then let $m_{s2} \leftarrow m_{s2}+1$
                      for any $s$ such that $ r \geq r_s \geq \ell (q)$.

              \end{itemize}

      \item  Let $\widehat{\Bbb{P}}_i  \leftarrow \frac{m_{i2}}{N}$
      and construct the confidence interval of confidence level $100(1- \delta) \%$.
\end{itemize}
\end{itemize}
\unskip

  It follows that $q$ can be viewed as
  an observation of a random variable
with uniform distribution over ${\cal B}(r_j)$ if and only if $r \geq r_j \geq \ell (q)$.
Hence, if the robustness requirement ${\bf P}$ has been
evaluated for ${\cal B}(r_i)$ at sample $q$, the result can be
accepted without repeated evaluation of {\bf P} for all ${\cal
B}(r_j)$ such that $ r \geq r_j \geq \ell (q)$.
{\it Thus, sample reuse allows us to save both the
sample generation and the evaluation of ${\bf P}$
for the sample.  It is also interesting to point out that
the samples collected for each $r_i$ are i.i.d.
and thus the confidence interval can be rigorously constructed based on
the evaluation of ${\bf P}$
for the samples.}

\section{Sample reuse factor}

Let ${\bf n}_i$ be the number of simulations required at $r_i$.
Define {\it sample reuse factor} as follows:\vspace{3pt}
\[
{\cal F}_{reuse} := \frac{ N l } { {\cal E} [ \sum_{i=1}^l  {\bf n}_i ] },
\]
\vspace{3pt}  where ${\cal E}(X)$ denotes the expectation of random variable $X$.
Obviously,  ${\cal F}_{reuse}$ measures the improvement of efficiency upon the conventional method.
We demonstrate that the improvement can be significant in most applications.

\begin{figure}[htbp]
\centerline{\psfig{figure=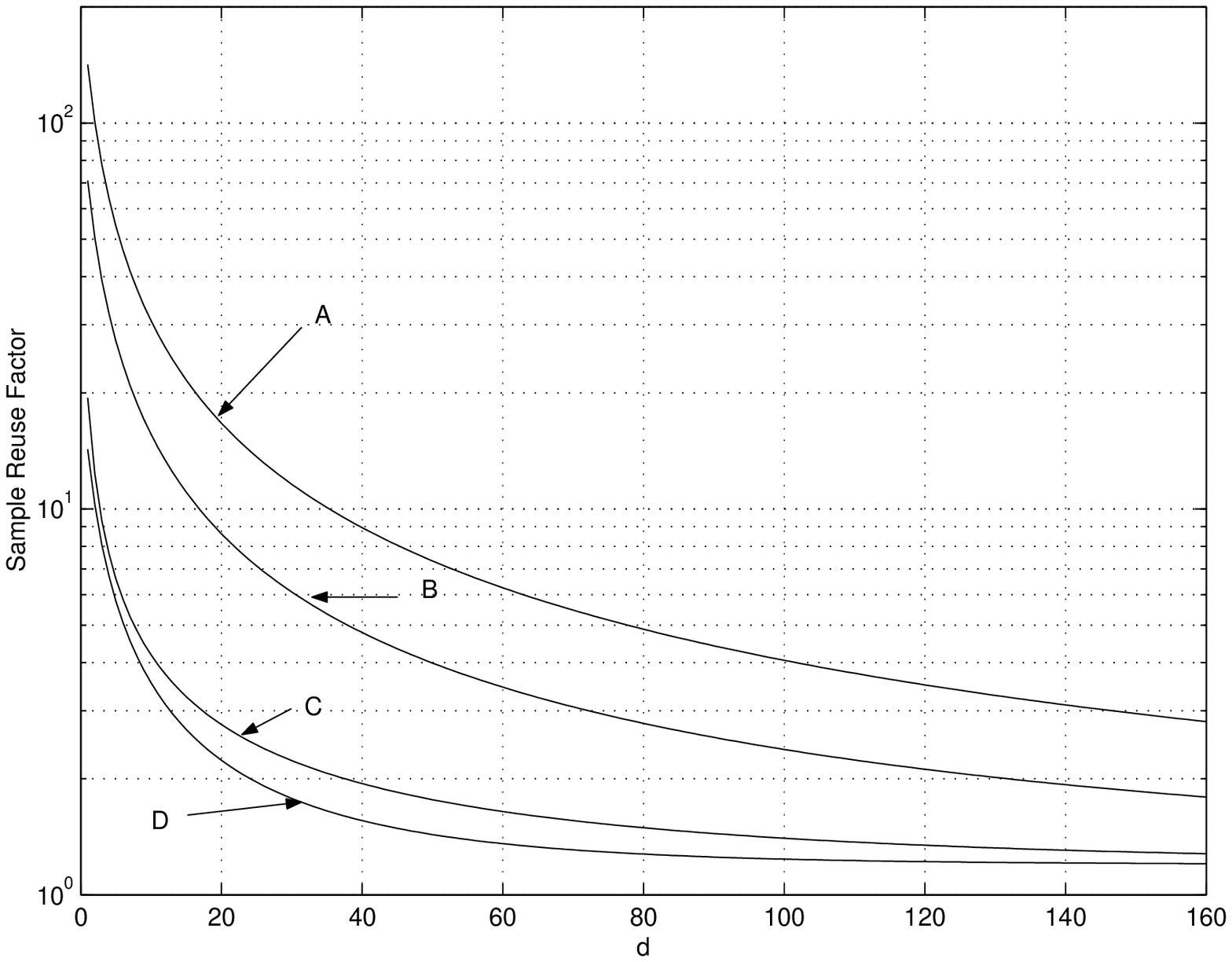, height=2in, width=4in
}} \caption{Performance improvement
(${\rm A:} \;l=200, \; b = 2 a; \;\;
{\rm B:} \;l=100, \; b = 2 a; \;\;
{\rm C:} \;l=100, \; a = 0; \;\;
{\rm D:} \;l=20, \; b = 2 a$).}
\label{fig_1}
\end{figure}

\begin{figure}[htbp]
\centerline{\psfig{figure=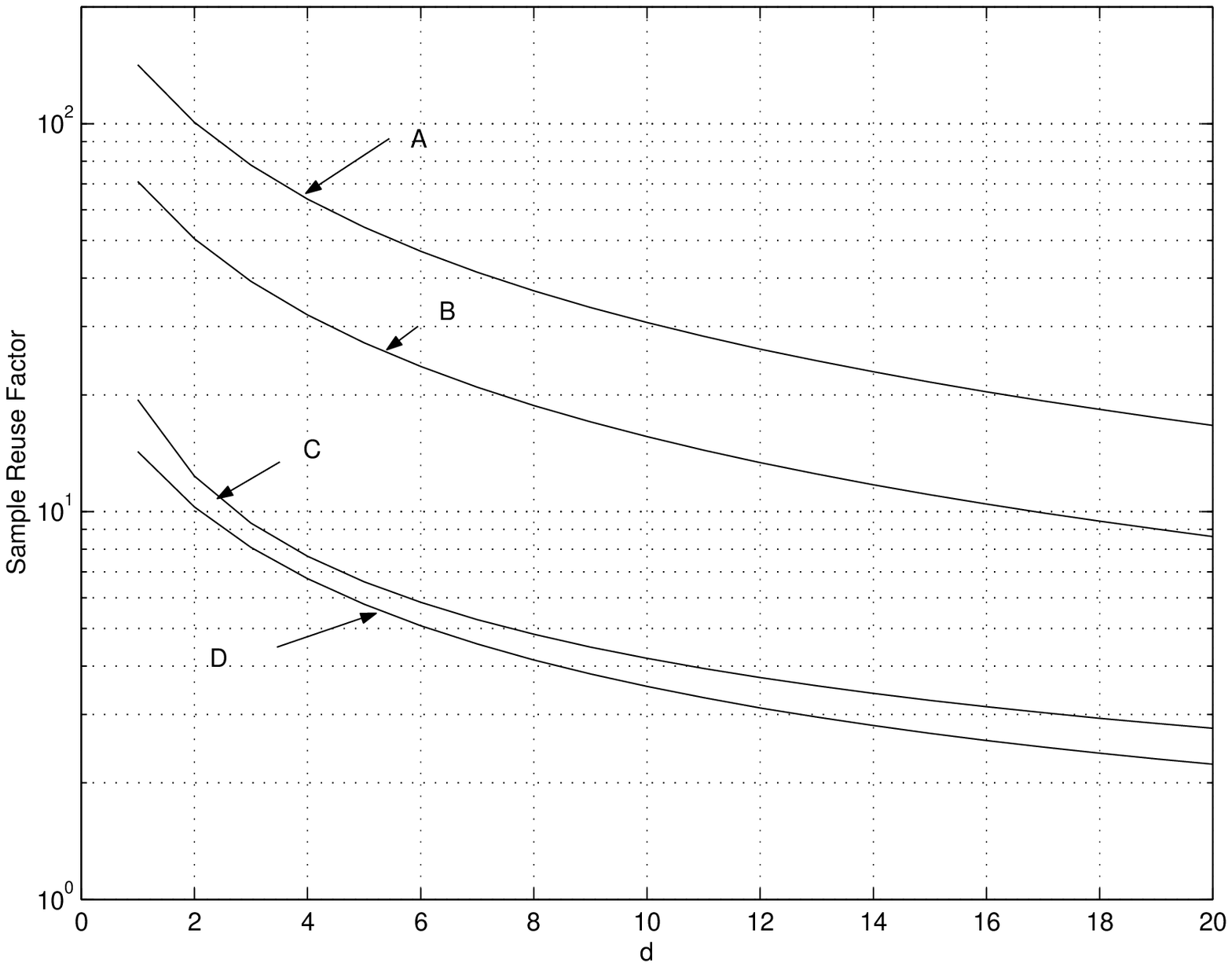, height=2in, width=4in
}} \caption{Performance improvement
(${\rm A:} \;l=200, \; b = 2 a; \;\;
{\rm B:} \;l=100, \; b = 2 a; \;\;
{\rm C:} \;l=100, \; a = 0; \;\;
{\rm D:} \;l=20, \; b = 2 a$).}
\label{fig_2}
\end{figure}

\begin{theorem} \label{reuse}
The {\it sample reuse factor}
${\cal F}_{reuse} = {l}/ {l -  \sum_{i=2}^l \left( \frac{r_i}{r_{i-1}} \right)^d}$,
where $d = n$ for $l_p$ ball
${\cal B}_p(r)$ and homogeneous star-shaped bounding set
${\cal B}_H (r)$; and \vspace{3pt}
\[
d = \sum_{i=1}^s \kappa(q_i) + \sum_{j=1}^c \kappa(\Delta_j)
\]
\vspace{3pt}  for spectral norm ball
${\cal B}_{\sigma}(r)$ with
$\kappa(.)$ defined as \vspace{3pt}
\[
\kappa(X) := \left\{\begin{array}{ll}
   2 m n  \;\;\;  {\rm if}\; X \; {\rm is \; a \;
   variable \; in} \; {\bf C}^{n \times m}\\
   m n \; \;\;\;
   {\rm if}\; X \; {\rm is \; a \; variable \; in} \; {\bf R}^{n \times m}.
\end{array} \right.
\]
\end{theorem}
\vspace{3pt} See the appendix for proof.  For illustration purposes,
we choose $r_i = b-\frac{(b-a)(i-1)}{l-1}$ for $i = 1,2, \ldots,l$.
By Theorem~\ref{reuse}, $ {\cal F}_{reuse} = {l}/ {l -  \sum_{i=2}^l
\left( 1 - \frac{1}{ \frac{l-1}{ 1 - \frac{a}{b} } - i + 2   }
\right)^d}$. Figures~\ref{fig_1} and ~\ref{fig_2} show that the
improvement over the conventional approach is significant when $d$
is not large. These figures also reveal that the sample reuse factor
does not scale well with the uncertainty dimension. For example,
when $d > 160$, the efficiency gained from sample reuse techniques
may not be attractive.

\section{Illustrative examples}
In this section we demonstrate through examples the power of the sample reuse algorithm
in solving a wide variety of complicated robustness analysis problems which are
intractable in the classical deterministic framework.

\begin{figure}[htbp]
\centerline{\psfig{figure=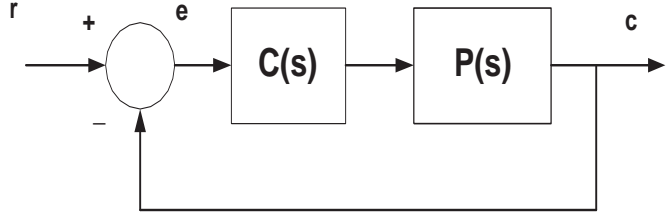, height=1.3in, width=4in
}} \caption{Uncertain system.}
\label{fig_08}
\end{figure}

First, we consider an example which has been studied in \cite{GS} by a deterministic approach.
The system is as shown in Figure~\ref{fig_08}.

The compensator is $C(s)=\frac{s+2}{s+10}$ and the plant is
$P(s)=\frac{800(1+0.1\delta_1)}{s(s+4+0.2\delta_2)(s+6+0.3\delta_3)}$
with parametric uncertainty $\Delta=[\delta_1,\delta_2,\delta_3]^{\rm T}$.
The nominal system is stable.
The closed-loop roots of the nominal system are \vspace{3pt}
\[
\hspace{-.40in}z_1= -15.9178,  \;\;  z_2 = -1.8309,  \;\; z_3 = -1.1256 + 7.3234i, \;\; z_4 = -1.1256 - 7.3234i.
\]
\vspace{3pt} The $H_\infty$
norm of the nominal closed-loop transfer function is ${||T^0||}_\infty = 2.78$.
The peak value, rise time, and settling time of step response of the
nominal system  are, respectively,
$P_{peak}^0 =1.47, \;\;  t_r^0 = 0.185, \;\mbox{{\rm and} } t_s^0 = 3.175$.
In all of the following examples, we take $l = 100$.
To guarantee that the absolute error of the estimate for the proportion
is less than $0.01$ with confidence level $99\%$, we choose $N = 26,492$
based on the well-known Chernoff bound (see \cite{KT,TD} for ``sharper'' bounds).
Since the Chernoff bound is conservative,
we also performed a post-experimental evaluation of the estimates by constructing confidence intervals
with confidence level $99\%$ based on Clopper--Pearson's method \cite{Clo}.

\begin{figure}[htbp]
\centerline{\psfig{figure=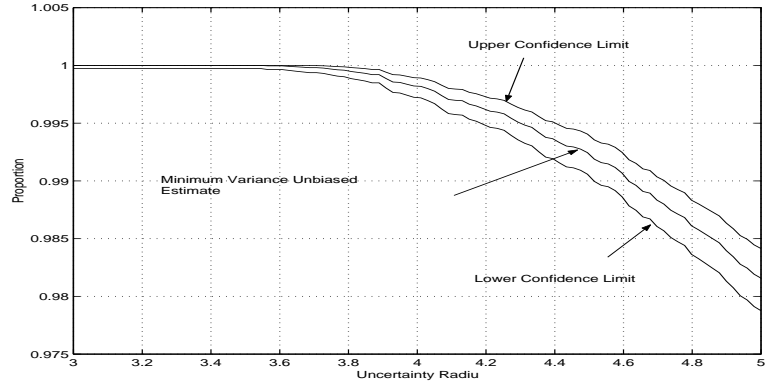, height=2in, width=4in
}} \caption{Robustness degradation curve
(reuse factor $= 41$). }
\label{fig_0}
\end{figure}

Figure~\ref{fig_0} is the robustness degradation curve
for robust stability over uncertainty set
${\cal B}_{\infty} (r) :=
\{ \Delta:\;||\Delta||_{\infty} \leq r \}$.
It demonstrates that a significant enhancement
of the robustness margin can be achieved at the price of a small risk.

\begin{figure}[htbp]
\centerline{\psfig{figure=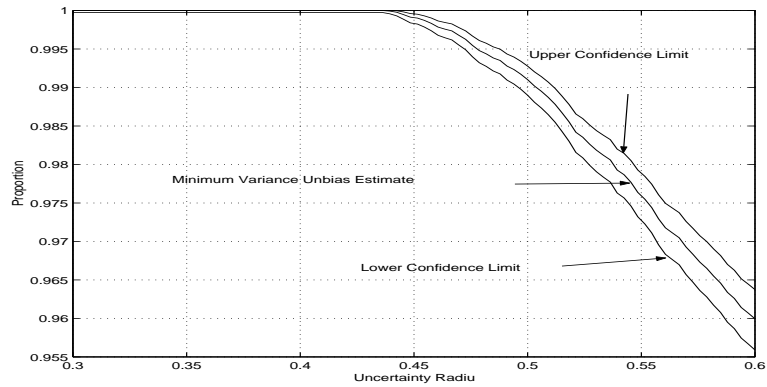, height=2in, width=4in
}} \caption{Robustness degradation curve
(reuse factor $= 43$). }
\label{fig_00}
\end{figure}

Figure~\ref{fig_00} is the robustness degradation curve, with the
robustness requirement ${\bf P}$ defined as stability and $H_\infty$
norm $< 170 \% \; {||T^0||}_\infty $, and the uncertainty set
defined as the ellipsoid ${\cal B}_2(r) :=\left\{
\Delta:\;||\Delta||_2 \leq r \right\}$.

\begin{figure}[htbp]
\centerline{\psfig{figure=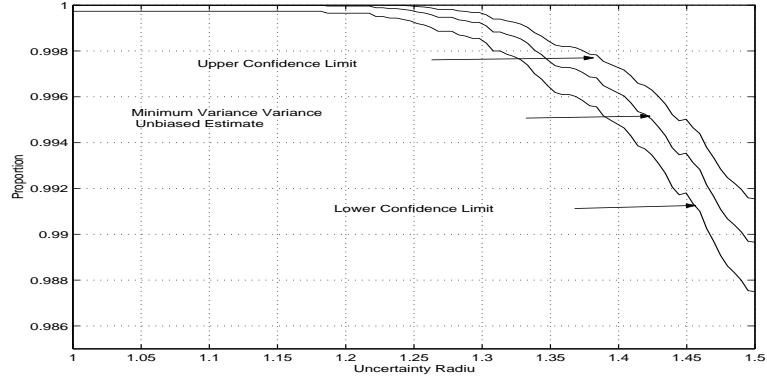, height=2in, width=4in
}} \caption{Robustness degradation curve
(reuse factor $= 49$).}
\label{fig_000}
\end{figure}

Figure~\ref{fig_000} is the robustness degradation curve with the
robustness requirement ${\bf P}$ defined as ${\cal D}$-stability
with the domain of poles defined as: real part $< -1.5$, or it falls
within one of the
 two disks centered at $z_3$
and $z_4$ with radius $0.3$.  The uncertainty set is defined as the polytope \vspace{3pt}
\[
{\cal B}_H (r):=\left\{ r \Delta +(1-r)\;\frac{ \sum_{i=1}^4 \Delta^i} {4} :
 \Delta \in {\rm conv}\{\Delta^1,\Delta^2,\Delta^3,\Delta^4\}  \right\},
\]
\vspace{3pt} where ``conv'' denotes the convex hull of $\Delta^i=[ \frac{1}{2} \sin(\frac{2i-1}{3}\pi),\;
\frac{1}{2} \cos(\frac{2i-1}{3}\pi),\;
-\frac{\sqrt{3}}{2}]^{\rm T}$ for $i=1,2,3$ and
$\Delta^4=[
0,\;0,\;1]^{\rm T}$.

\begin{figure}[htbp]
\centerline{\psfig{figure=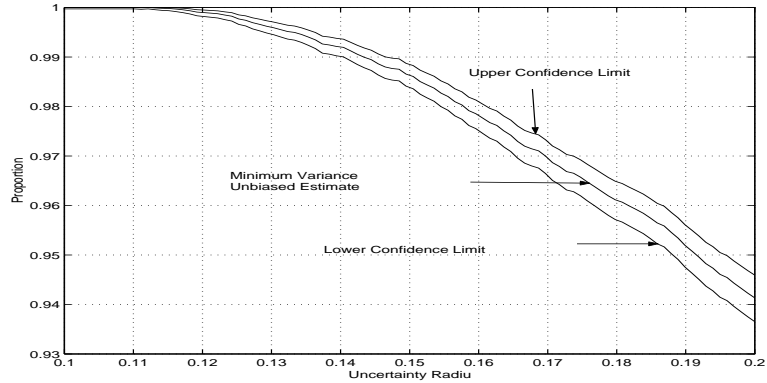, height=2in, width=4in
}} \caption{Robustness degradation curve
(reuse factor $= 38$).}
\label{fig_0000}
\end{figure}

Figure~\ref{fig_0000} is the robustness degradation curve for the
case where the uncertainty set is ${\cal B}_{\infty} (r) := \{
\Delta:\;||\Delta||_{\infty} \leq r \}$, the robustness requirement
${\bf P}$ is: stability, and rise time $t_r < 135 \% \; t_r^0  =
0.25$, settling time $t_s < 110 \% \; t_s^0 = 3.5$, and overshoot
$P_{peak} < 116 \% \; P_{peak}^0= 1.7$.

\begin{figure}[htbp]
\centerline{\psfig{figure=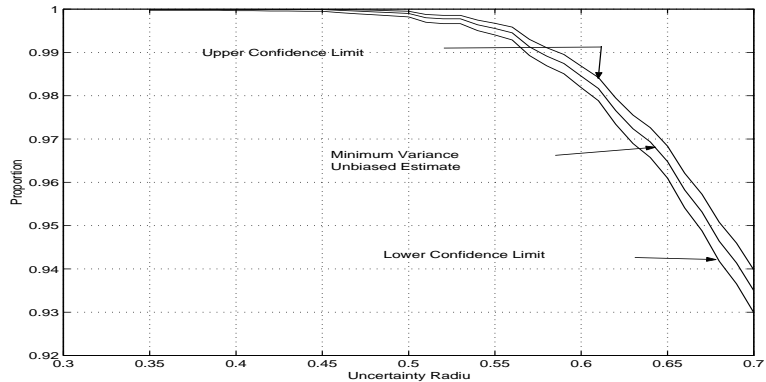, height=2in, width=4in
}} \caption{Robustness degradation curve
(reuse factor $= 5$).}
\label{fig_00000}
\end{figure}

Finally, we consider the same example in \cite{CDT1} where the class of uncertainty is defined as \vspace{3pt}
\[
{\bf \Delta} := \{ {\rm blockdiag} [q_1 I_5, \; q_2 I_5,\; \Delta_1] \},
\]
\vspace{3pt}  where $\Delta_1 \in {\bf C}^{4 \times 4}$, and $I_5$ denotes the identity matrix of $5 \times 5$.
By Theorem~\ref{reuse}, we have $d = 34$.
Figure~\ref{fig_00000} shows the robustness degradation curve.
An improvement (of efficiency) about fivefold is achieved by our algorithm.

\section{A probabilistic perspective}

In sections 2 and 3, we have derived our sample reuse
algorithm from the worst-case deterministic framework.
In this section, we show that the proposed algorithm is also applicable from the
perspective of the random nature of uncertainty.
In situations where we need to take into account the random nature of uncertainty,
the pioneering work of Barmish, Lagoa, Tempo, Bai, and Fu \cite{BLT, bai} allows our
sample reuse algorithm to be applied to
solve efficiently a wide variety of robustness problems.
The following theorem plays an important role.

\begin{theorem}  [{\rm see \cite {BLT}}] \label{bar}
Suppose that uncertainty ${\bf \Delta}$ is a random variable with a density function
$f(\Delta)$ which depends only on $\ell ( \Delta)$ and is nonincreasing
with respect to $\ell (\Delta)$.  Then \vspace{3pt}
\[
\Pr \{
{\rm The \; associated \; system}
\; G({\bf \Delta}) \; {\rm guarantees} \; {\bf P}  \;\; | \;\; {\bf \Delta} \in {\cal B}( r) \}
 \geq \inf_{0 \leq \rho \leq r} \Bbb{P} (\rho).
\]
\end{theorem}
\vspace{3pt}
{\it Remark {\rm 1}.} A remarkable fact of Theorem~\ref{bar} is that no assumption
 needs to be imposed on the robustness requirement ${\bf P}$.
 The assumption in Theorem~\ref{bar} is roughly interpreted to mean that
the probability measure of the uncertainty is radially symmetrical with respect to the nominal value.
In many applications, small perturbations are
 more likely than large perturbations, and
 the uncertainty is sufficiently unstructured so as to be treated
 equally likely in the surface of ${\cal B}(r)$ \cite{BLT}.


 {\it Remark {\rm 2}.} It should be noted that
  Theorem~\ref{bar} applies to a homogeneous star-shaped bounding set, $l_p$ ball, and spectral norm ball.
 We introduce in Theorem~\ref{bar} a conditional probability based on the following reason:  It does not seem logical
 to treat the uncertainty as different bounded random variables.
 For example, if the uncertainty possesses a certain distribution over ${\cal B}( r_1)$,
it would be a contradiction that the uncertainty possesses another
distribution over ${\cal B}( r_2)$ for $r_2 > r_1$. In fact, if the
uncertainty is of random nature, then the associated distribution is
unique.

Based on Theorem~\ref{bar}, we can apply the sample reuse algorithm to
estimate $\Bbb{P}(r)$ for $r \in [0, b]$, from which we can construct the lower bounds for
$\Pr \{ G({\bf \Delta}) \; {\rm guarantees} \\  \; {\bf P}  \;\; | \;\; {\bf \Delta} \in {\cal B}( r) \}$.

\section{Conclusion}

We develop a fast algorithm for computing the robustness degradation function
which overcomes the computational complexity and
conservatism issue of the deterministic worst-case methods.  We also demonstrate that our algorithm
can provide efficient solutions for a wide variety of robustness analysis problems
which are intractable by the deterministic  worst-case methods.  We derive our algorithm from
the worst-case deterministic framework and also show that
the algorithm is applicable from a probabilistic perspective.

\appendix

\section{Proof of Theorem 1}

The following lemma follows essentially from the definition of volume function ${\rm vol} (. )$.
\begin{lemma}   \label{lem1}
Let $X = \{ \Delta \in {\bf C}^{n \times m}: \bar{\sigma}(\Delta) \leq r \}$ and
$Y = \{ \Delta \in {\bf C}^{m \times n}: \bar{\sigma}(\Delta) \leq r \}$.
Let $Z = \{ \Delta \in {\bf R}^{n \times m}: \bar{\sigma}(\Delta) \leq r \}$ and
$W = \{ \Delta \in {\bf R}^{m \times n}: \bar{\sigma}(\Delta) \leq r \}$.
Then ${\rm vol}(X) = {\rm vol}(Y)$ and ${\rm vol}(Z) = {\rm vol}(W)$.
\end{lemma}

\begin{lemma} \label{lem2}
Let $m \geq n$.  Define
spectral norm ball ${\cal B}_{\sigma}^C(r) = \{ \Delta \in {\bf C}^{n \times m}: \bar{\sigma}(\Delta) \leq r \}$
and spectral norm ball
${\cal B}_{\sigma}^R(r) = \{ \Delta \in {\bf R}^{n \times m}: \bar{\sigma}(\Delta) \leq r \}$.
Then
${\rm vol} ( {\cal B}_{\sigma}^C(r) )  = {\rm vol} ( {\cal B}_{\sigma}^C(1) ) r^d$ with
$d = 2 m n$  and
${\rm vol} ( {\cal B}_{\sigma}^R(r) )  = {\rm vol} ( {\cal B}_{\sigma}^R(1) ) r^d$ with
$d = m n$.

\end{lemma}

\bpf
 By Theorem 1 of \cite{CDT1}, we have \vspace{3pt}
\[ \frac{ \Upsilon_C } { {\rm vol} ( {\cal B}_\sigma^C(r)   ) } \;
\int_{r \geq \sigma_1 > \sigma_2 > \cdots > \sigma_n > 0 }
\prod_{i=1}^n \sigma_i^{2(m-n)+1} \times \prod_{1 \leq i < k \leq n} (\sigma_i^2 - \sigma_k^2)^2  \; d \sigma_1 d \sigma_2 \cdots d \sigma_n = 1\\
\]
\vspace{3pt}with $\Upsilon_C= \frac{2^n \pi^{mn}   } {  \prod_{k=1}^n (n-k)! (m-k)! }$.
Performing a change of variables with $x_i = \frac{\sigma_i}{r}$ for $i = 1, \ldots,n$, we have \vspace{3pt}
\[
\frac{ r^{2 m n}\Upsilon_C } { {\rm vol} ( {\cal B}_\sigma^C(r)   )
} \; \int_{1 \geq x_1 > x_2 > \cdots > x_n > 0   }
\prod_{i=1}^n x_i^{2(m-n)+1} \\
\times \prod_{1 \leq i < k \leq n} (x_i^2 - x_k^2)^2  \; d x_1 d x_2
\cdots d x_n = 1.
\]
Thus \[ {\rm vol} ( {\cal B}_\sigma^C(1)   ) = \Upsilon_C \int_{1
\geq x_1 > x_2 > \cdots > x_n > 0   } \prod_{i=1}^n x_i^{2(m-n)+1}
\times \prod_{1 \leq i < k \leq n} (x_i^2 - x_k^2)^2  \; d x_1 d x_2
\cdots d x_n,
\]
and
\[
{\rm vol} ( {\cal B}_\sigma^C(r)   ) = {\rm vol} ( {\cal B}(1)   ) \; r^{2 m n}.
\]

Similarly, by Theorem 2 of \cite{CDT1}, we can show that ${\rm vol}
( {\cal B}_{\sigma}^R(r) )  = {\rm vol} ( {\cal B}_{\sigma}^R(1) )
r^d$ with $d = m n$. \qquad \epf

\begin{lemma}  \label{lem3}
${\rm vol} ( {\cal B}(r) ) = {\rm vol} ( {\cal B}(1) ) r^d$ where $d = n$ for $l_p$ ball
${\cal B}_p(r)$ and homogeneous star-shaped bounding set ${\cal B}_H (r)$; and
$d = \sum_{i=1}^s \kappa(q_i) + \sum_{j=1}^c \kappa(\Delta_j)$
for spectral norm ball ${\cal B}_{\sigma}(r)$.
\end{lemma}

\bpf The truth is obvious for cases of an $l_p$ ball and homogeneous
star-shaped bounding set. To prove the lemma for the case of a
spectral norm ball, we need to apply Lemmas~\ref{lem1} and
\ref{lem2}. \qquad \epf

\begin{lemma} \label{lem4}
For $i = 1, \ldots, l-1$, \vspace{3pt}
\[
{\cal E}\left[ {\bf n}_{i+1} \right] =
N - \sum_{j=1}^i \left( \frac{r_{i+1}}{r_{j}} \right)^d \; {\cal E}\left[ {\bf n}_{i} \right] .
\]
\end{lemma}
\vspace{3pt}

\bpf

 Let $q^1, q^2, \ldots, q^{ {\bf n}_j }$ be the samples
generated on $r_j$. For $l = 1, \ldots, {\bf n}_j$, define binomial
random variable $X_{j, i+1}^l$ such that \vspace{3pt}
\[
X_{j, i+1}^l := \left\{\begin{array}{ll}
   1 \;\;\;  {\rm if}\; q^l \; {\rm fall \; in} \; {\cal B}(r_{i+1}),\\
   0 \; \;\;\;
   {\rm otherwise.}
\end{array} \right.
\]
\vspace{3pt}  By the rule of the sample reuse algorithm, \vspace{3pt}
\[
N = {\bf n}_{i+1} + \sum_{j=1}^i \sum_{l=1}^{ {\bf n}_j  } X_{j, i+1}^l.
\]
Thus for $i = 1, \ldots, l-1$,
\[
{\cal E}\left[ {\bf n}_{i+1} \right]  = N - \sum_{j=1}^i {\cal E}
\left[ \sum_{l=1}^{ {\bf n}_j  } X_{j, i+1}^l \right] =  N -
\sum_{j=1}^i  \sum_{n \in \Omega_{ {\bf n}_j  }} \sum_{l=1}^n {\cal
E} \left[ X_{j, i+1}^l  \; |  \; {\bf n}_j  = n \right] \; \Pr \{
{\bf n}_j = n  \},
\]
where $\Omega_{ {\bf n}_j  }$ denotes the sample space of ${\bf
n}_j$. Since $q^l$ is a random variable with uniform distribution
over ${\cal B}(r_i)$, it follows from Lemma~\ref{lem3} that
\vspace{3pt}
\[
{\cal E} [ X_{j, i+1}^l  \; |  \; {\bf n}_j  = n ] =
\frac{ {\rm vol} ({\cal B}(r_{i+1}))  }
{ {\rm vol} ({\cal B}(r_j) ) } = \left( \frac{r_{i+1}}{r_{j}} \right)^d.
\]
Therefore, \bee {\cal E}\left[ {\bf n}_{i+1} \right]
 & = & N - \sum_{j=1}^i  \sum_{n \in \Omega_{ {\bf n}_j  }}  n \; \left( \frac{r_{i+1}}{r_{j}} \right)^d \;
\Pr \{ {\bf n}_j = n  \}\\
& = & N - \sum_{j=1}^i  \; \left( \frac{r_{i+1}}{r_{j}} \right)^d \;
\sum_{n \in \Omega_{ {\bf n}_j  }}  n \Pr \{ {\bf n}_j = n  \}\\
 & = & N - \sum_{j=1}^i \left( \frac{r_{i+1}}{r_{j}} \right)^d \; {\cal
E}\left[ {\bf n}_{i} \right]. \eee

\epf

\vspace{3pt}

\begin{lemma} \label{lem5}
For $i = 2, \cdots, l$, \vspace{3pt}
\[
{\cal E}\left[ {\bf n}_i \right] = N - N  \left( \frac{r_i}{r_{i-1}} \right)^d.
\]
\end{lemma}
\vspace{3pt} \bpf We use induction.  Obviously, \vspace{3pt}
\[
{\cal E}\left[ {\bf n}_{1} \right] = N.
\]
\vspace{3pt} By Lemma~\ref{lem4}, we get \vspace{3pt}
\[
{\cal E}\left[ {\bf n}_{2} \right] = N - N \left( \frac{r_{2}}{r_{1}} \right)^d.
\]
\vspace{3pt} Suppose it is true that \vspace{3pt}
\[
{\cal E}\left[ {\bf n}_i \right] = N - N  \left( \frac{r_i}{r_{i-1}} \right)^d.
\]
\vspace{3pt} Then \vspace{3pt}
\begin{eqnarray*}
\sum_{j = 1}^i \; \left( \frac{r_{i+1}}{r_{j}} \right)^d \; {\cal
E}\left[ {\bf n}_j \right] & =  &\sum_{j = 1}^i \left(
\frac{r_{i+1}}{r_{j}} \right)^d\; \left[ N - N  \left(
\frac{r_j}{r_{j-1}} \right)^d  \right] \\ \vspace{3pt} & = & \sum_{j
= 1}^i \left[ N \left( \frac{r_{i+1}}{r_{j}} \right)^d
- N  \left( \frac{r_{i+1}}{r_{j-1}} \right)^d  \right]\\
& = & N \left( \frac{r_{i+1}}{r_{i}} \right)^d.
\end{eqnarray*}
\vspace{3pt} It follows from Lemma~\ref{lem4} that \vspace{3pt}
\[
{\cal E}\left[ {\bf n}_{i+1} \right]
 =  N - \sum_{j = 1}^i \; \left( \frac{r_{i+1}}{r_{j}} \right)^d \; {\cal E}\left[ {\bf n}_j \right]
 =  N - N \left( \frac{r_{i+1}}{r_{i}} \right)^d.
\]
\vspace{3pt} The proof of Lemma~\ref{lem5} is thus completed by
induction. \qquad \epf

Now we are in the position to prove Theorem~\ref{reuse}.  By Lemmas~\ref{lem4} and \ref{lem5},
we have \vspace{3pt}
\[
{\cal E} \left[ \sum_{i=1}^l  {\bf n}_i \right] =
N + \sum_{i=2}^l \left[ N - N  \left( \frac{r_i}{r_{i-1}} \right)^d \right]
= N l - N \; \sum_{i=2}^l \left[ 1 - \left( \frac{r_i}{r_{i-1}} \right)^d \right].
\]
\vspace{3pt} Therefore, \vspace{3pt}
\[
{\cal F}_{reuse}  =  \frac{ N l } { {\cal E} [ \sum_{i=1}^l  {\bf
n}_i ] }
 =  \frac{l} {l -  \sum_{i=2}^l \left( \frac{r_i}{r_{i-1}} \right)^d}
\]
\vspace{3pt} and thus the proof of Theorem~\ref{reuse} is completed.
\qquad  \pagebreak

\end{document}